\documentclass[11pt]{article}

\usepackage{amssymb}
\usepackage{amsmath}
\usepackage{amsthm}
\usepackage{enumerate}

\newtheoremstyle{plain}%
    {8pt plus2pt minus4pt}%
    {8pt plus2pt minus4pt}%
    {\itshape}%
    {}%
    {\bfseries\scshape}%
    {}%
    {6pt}
    {}%

\newtheoremstyle{remark}%
    {8pt plus2pt minus4pt}%
    {8pt plus2pt minus4pt}%
    {\upshape}
    {}%
    {\bfseries\scshape}%
    {}%
    {6pt}
    {}%

\theoremstyle{plain}
\newtheorem{thm}{Theorem}[section]

\theoremstyle{remark}

\title{On generalized Ramsey numbers for 3-uniform hypergraphs}

\author{
{\Large{Andrzej Dudek}}\thanks{
\footnotesize {Department of Mathematics, Western~Michigan~University, Kalamazoo, MI 49008, \texttt{andrzej.dudek@wmich.edu}
}}
\and
{\Large{Dhruv Mubayi}}\thanks{
\footnotesize {Department of Mathematics, Statistics, and Computer Science, University of Illinois at Chicago, Chicago, IL 60607, \texttt{mubayi@math.uic.edu }}}~\thanks{\footnotesize{Research partially supported by NSF grant DMS-0969092}}
}

\date{\today}

\begin{document}
\maketitle		
		
\begin{abstract}
The well-known Ramsey number $r(t,u)$ is the smallest integer~$n$ such that every $K_t$-free graph of order~$n$ contains an independent set of size~$u$. In other words, it contains a subset of $u$~vertices with no~$K_2$. Erd{\H o}s and Rogers introduced a more general problem replacing~$K_2$ by~$K_s$ for $2\le s<t$. Extending the problem of determining Ramsey numbers they defined the numbers
$$
f_{s,t}(n)=\min \big{\{} \max \{ |W| : W\subseteq V(G) \text{ and } G[W] \text{ contains no } K_s\}\big{\}},
$$
where the minimum is taken over all $K_t$-free graphs~$G$ of order~$n$. In this note, we study an analogous function $f_{s,t}^{(3)}(n)$ for 3-uniform hypergraphs. In particular, we show that there are constants $c_1$ and $c_2$ depending only on $s$ such that

$$
c_1(\log n)^{1/4} \left(\frac{\log\log n}{\log\log\log n}\right)^{1/2} < f_{s, s+1}^{(3)}(n) < c_2 \log n.
$$

\end{abstract}

\section{Introduction}
A {\em hypergraph} $\mathcal{G}$ is a pair $(V,\mathcal{E})$ such that $V = V(\mathcal{G})$ is a set of vertices and $\mathcal{E}\subseteq 2^V$ is a set of hyperedges.
A {\em $k$-uniform hypergraph} (also called a {\em $k$-graph}) is a hypergraph such that all its hyperedges have size $k$. Denote by $K_{t}^{(k)}$ the complete $k$-uniform hypergraph of order~$t$.
Let the {\em Ramsey number} $r^{(k)}(t,u)$ denote the smallest integer $n$ such that any red-blue coloring of the edges of $K_{n}^{(k)}$ yields a red copy of $K_{t}^{(k)}$ or a blue copy of $K_{u}^{(k)}$. It is well-known due to a result of Ramsey~\cite{RA} that such numbers are finite. In other words, $r^{(k)}(t,u)$ is the smallest integer~$n$ such that every $K_t^{(k)}$-free hypergraph of order~$n$ contains an independent set of size~$u$, or equivalently, it contains a $u$-subset of vertices with no~$K_k^{(k)}$. One can consider a more general problem replacing~$K_k^{(k)}$ by~$K_s^{(k)}$ for some $k\le s<t$.
For fixed integers $k \le s < t$ let
$$
f_{s,t}^{(k)}(n)=\min \big{\{} \max \{ |W| : W\subseteq V(\mathcal{G}) \mbox{ and } \mathcal{G}[W] \mbox{ contains no } K_s^{(k)}\} \big{\}},
$$
where the minimum is taken over all $K_{t}^{(k)}$-free $k$-uniform hypergraphs~$\mathcal{G}$ of order~$n$. To prove the lower bound $f_{s,t}^{(k)}(n)\ge a$ one must show that every $K_{t}^{(k)}$-free $k$-graph of order~$n$ contains a subset of $a$ vertices with no copy of~$K_s^{(k)}$. To prove the upper bound $f_{s,t}^{(k)}(n) < b$ one must construct a $K_{t}^{(k)}$-free $k$-graph of order~$n$  such that every subset of $b$ vertices contains a copy of~$K_s^{(k)}$.

As we already observed the problem of determining $f_{s,t}^{(k)}(n)$ extends that of finding Ramsey numbers. Formally,
$$
r^{(k)}(t,u) = \min \{ n : f_{k,t}^{(k)}(n) \ge u\}.
$$

For graphs (i.e. $k=2$) the above function was first considered by Erd{\H o}s and Rogers~\cite{ERog} only for $t=s+1$, which might be viewed as the most restrictive case. Since then the function has been studied by several researchers including Alon, Bollob\'as, Dudek, Erd\H{o}s, Gallai, Hind, Krivelevich, Retter, R\"odl, Sudakov, and Wolfovitz \cite{AK, BH, DRR, DR1, DR2, EG, KR2, KR, SU, SU2, W}.  All logs in this paper are to the base $e$.
For any $s\ge 3$ the best published bounds are of the form
\begin{equation}\label{eq:graphs}
\Omega\big{(}(n \log \log n)^{1/2}\big{)} = f_{s, s+1}^{(2)}(n) = O\big{(} (\log n)^{4s^2} n^{1/2}\big{)},
\end{equation}
where the lower bound comes from~\cite{KR} and the upper bound from~\cite{DRR}.
However, using the  result due to Shearer~\cite{Shearer95}, that an $n$-vertex  $K_{s}$-free graph with average degree $d$ has an independent set of size at least
$$
\Omega \left( \frac{n}{d} \frac{\log d}{\log\log d}\right),
$$
we can immediately improve the lower bound in (\ref{eq:graphs}) to a new bound
\begin{equation} \label{new} 
\Omega\left(\left(\frac{n \log n }{\log\log n}\right)^{1/2}\right).
\end{equation}
Indeed, if there is a vertex of degree at least $(\frac{n \log n }{\log\log n})^{1/2}$ in a $K_{s+1}$-free graph, then its neighborhood is $K_{s}$-free and we are done. Otherwise, the average degree is less than $(\frac{n \log n }{\log\log n})^{1/2}$ and by the Shearer result there is an independent set of size bounded from below by \eqref{new}.

In this note we extend the previous results to  3-uniform hypergraphs and show the following bounds.
\begin{thm}\label{thm:3}
For all integers $3\le s < t$ and every $n$
$$
f_{s-1,t-1}^{(2)}\big{(}\lfloor \sqrt{\log n} \rfloor\big{)}  \le f_{s,t}^{(3)} (n) \le C \log n,
$$
where  $C$ is a positive constant depending only on $s$.
\end{thm}
\noindent
In particular, for $t=s+1$ the above theorem together with \eqref{new} implies that there are positive constants $c_1$ and $c_2$ depending only on $s$ such that
$$
c_1(\log n)^{1/4} \left(\frac{\log\log n}{\log\log\log n}\right)^{1/2} < f_{s, s+1}^{(3)}(n) < c_2 \log n.
$$

\section{General upper bound}

For $s\ge k\ge 3$, we show that there is a $K_{s+1}^{(k)}$-free $k$-graph $\mathcal{G}$ such that any subset of its vertices of size $C(\log n)^{1/(k-2)}$ (for some $C=C(k,s)$) contains $K_{s}^{(k)}$. This will imply that
\begin{equation}\label{eq:k}
f_{s,s+1}^{(k)} (n) \le C  (\log n)^{1/(k-2)}.
\end{equation}
Moreover, since for any $t\ge s+1$
$$
f_{s,t}^{(k)} (n) \le f_{s,s+1}^{(k)} (n),
$$
setting $k=3$ will yield the upper bound in Theorem~\ref{thm:3}.

The following construction has some similarities to the approach taken by Conlon, Fox and Sudakov in the proof of Theorem~3.1 in~\cite{CFS}.
Let $\chi : \binom{[n]}{k-1} \to [s-k+2]$ be a coloring of all $(k-1)$-tuples chosen uniformly at random from all $(s-k+2)$-colorings of $\binom{[n]}{k-1}$. That means every $(k-1)$-tuple is colored with one particular color with probability $1/(s-k+2)$ independently of all other $(k-1)$-tuples. Now we construct $\mathcal{G} = \mathcal{G}(\chi)$. Let $V=\{1,2,\dots,n\}$ be the vertex set of $\mathcal{G}$. The set of hyperedges consists of all $k$-tuples $e = \{u_1<\dots<u_{k-1}<u_k\}$ for which $\chi(e\setminus u_k) \neq \chi(e \setminus u_{k-1})$.

First observe that $\mathcal{G}$ is $K_{s+1}^{(k)}$-free. Indeed, let $\{v_1<v_2<\dots<v_{s+1}\}$ be a set of vertices inducing a clique $K_{s+1}^{(k)}$. Then, in particular, all $k$-tuples $\{v_1<\dots<v_{k-1}<v_{i}\}$ for $k\le i\le s+1$ must be present. Hence, all $s-k+3$ of the numbers $\chi(v_1,\dots, v_{k-2},v_{i})$, $k-1\le i\le s+1$, must be pairwise different. Clearly, this gives a contradiction, since we only have $s-k+2$ colors.

Now note that with a positive probability, say $p$, which does not depend on $n$, any fixed set $S=\{v_1<v_2<\dots<v_{s}\}$ of vertices will induce  $K_s^{(k)}$. Indeed, it suffices to observe that if for each $\{v_{j_1}<\dots<v_{j_{k-1}}\} \in \binom{S}{k-1}$
$$
\chi(v_{j_1},\dots,v_{j_{k-1}}) = j_{k-1}-k+2,
$$
then every $k$-tuple of $S$ is present.

Let $C(k-1,s,W)$ be a maximum collection of $s$-sets of set~$W$ such that no $(k-1)$-set in $W$ is covered more than once. Denote its size by
$c(k-1,s,W)$. Due to a result of R\"odl~\cite{Rodl} (see also Section 4.7 in~\cite{AS}) it is known that
$$
\lim_{|W|\to \infty} \frac{c(k-1,s,W)}{\binom{|W|}{k-1} / \binom{s}{k-1}} = 1,
$$
hence, $c(k-1,s,W) = \Theta(|W|^{k-1})$.

For a fixed subset $S \subseteq V$ of size $s$ let $A_S$ be the event that $\mathcal{G}[S]$ is a clique of size~$s$. As we already observed
$$
\Pr(A_S) = p
$$
for some $0<p<1$.
Let $W$ of size $w$ be a subset of vertices of $V$.
Thus,
$$
\Pr\left(\bigcap_{S\subseteq W} \bar{A}_S\right) \le \Pr\left(\bigcap_{S \in C(k-1,s,W)} \bar{A}_S\right).
$$
Moreover, since for all different  $S$ and $S'$ in  $C(k-1,s,W)$, $|S\cap S'| \le k-2$, and consequently, all events $A_S$ for $S\in C(k-1,s,W)$ are mutually independent. Thus,
$$
\Pr\left(\bigcap_{S\subseteq W} \bar{A}_S\right) \le \prod_{S \in C(k-1,s,W)} \Pr(\bar{A}_S) = (1-p)^{\Theta(w^{k-1})}
$$
and hence, the union bound yields
$$
\Pr\left(\bigcup_{W\subseteq V} \bigcap_{S\subseteq W} \bar{A}_S\right) \le \binom{n}{w} (1-p)^{\Theta(w^{k-1})} < 1
$$
provided that $w=\Omega( (\log n)^{1/(k-2)} )$. Thus, there is a coloring $\chi$ and a constant $C=C(k,s)$ such that every subset of $C(\log n)^{1/(k-2)}$ of vertices in $\mathcal{G}(\chi)$ contains a clique of size $s$, as required.

\section{Lower bound for $3$-graphs}\label{sec:lower}

In this section we prove the lower bound from Theorem~\ref{thm:3}.

We show that every $K_{t}^{(3)}$-free $3$-graph $\mathcal{G}$ of order~$n$ contains a subset of $f_{s-1,t-1}^{(2)}(\lfloor \sqrt{\log n} \rfloor)$ vertices  with no copy of~$K_s^{(3)}$. In order to do it, we will adapt the classical approach of Erd\H{o}s and Rado~\cite{ErdRado} (see also Section 1.2 in \cite{GRS}) for finding an upper bound on Ramsey numbers.

Given $n$, choose $m$ such that 
\begin{equation}\label{ineq:m}
2^{m^2} \le n < 2^{m^2+1}.
\end{equation} 
Let $\mathcal{G}=(V,\mathcal{E})$ be a  $K_{t}^{(3)}$-free $3$-graph of order $n$. Denote the complement of $\mathcal{E}$ (i.e. the set of all triples which are not in $\mathcal{E}$) by $\mathcal{E}^c$. We will greedily construct a sequence $A = \{v_1,\dots,v_{m},v_{m+1}\} \subseteq V$ such that for any given pair $1\le i<j\le m$ all triples $\{v_i,v_j,v_k\}$ with $j<k\le m+1$ are in $\mathcal{E}$ or all of them are in $\mathcal{E}^c$. Assume for a while that we can construct such a sequence~$A$. Let $G$ be a graph on set $\{v_1,\dots,v_m\}$ such that $\{v_i,v_j\}$, $i<j$,  is an edge if and only if all triples  $\{v_i,v_j,v_k\}$ with $j<k\le m+1$ are in $\mathcal{E}$. First observe that $G$ is $K_{t-1}$-free graph. Otherwise, if $S$ is a set of vertices of $G$ that induces $K_{t-1}$, then $S\cup \{v_{m+1}\}$ induces a clique $K_{t}^{(3)}$ in $\mathcal{G}$, a contradiction.

Since $G$ is $K_{t-1}$-free, we can find a subset $W$ of at least $f_{s-1,t-1}^{(2)}(m)$ vertices with no $K_{s-1}$. Now it is not difficult to see that the subhypergraph induced by $W$ in $\mathcal{G}$ contains no $K_{s}^{(3)}$. Otherwise, let $S\subseteq W$ be a set of size $s$ such that $\mathcal{G}[S]$ is a clique $K_{s}^{(3)}$. Moreover, let $v \in S$ be a vertex which appears latest in sequence~$A$ compared to other vertices in $S$. Then $S\setminus v \subseteq W$ induces a clique $K_{s-1}$ in $G$, a contradiction.

We just showed that every $K_{t}^{(3)}$-free graph $\mathcal{G}$ of order $n$ satisfying~\eqref{ineq:m} has a subhypergraph of order at least
$$
f_{s-1,t-1}^{(2)}(m)$$
 with no $K_s^{(3)}$. Since $n < 2^{m^2+1}$, we see that $m > \sqrt{\log n}$ as desired.

It remains to show how to construct sequence $A$, and this is the argument of Erd\H{o}s and Rado. Pick any vertex $v_1$ of $V=V(\mathcal{G})$. Set $V_1 = V\setminus \{v_1\}$. Assume that we already picked $\{v_1,\dots,v_i\}$ and set $V_i$ such that all triples $\{v_a,v_b,w\}$ with $1\le a<b\le i$ and $w\in V_i$ are in $\mathcal{E}$ or all are in $\mathcal{E}^c$. Choose any vertex from $V_i$, say $v_{i+1}$. Now we show how to define $V_{i+1}$ such that $|V_{i+1}| \ge (|V_i|-1) / 2^{i}$ and all triples $\{v_a,v_b,w\}$ with $1\le a<b\le i+1$ and $w\in V_{i+1}$ are in $\mathcal{E}$ or all are in $\mathcal{E}^c$.
Let $V_{i,0} = V_i \setminus \{v_{i+1}\}$. Suppose we already constructed $V_{i,j} \subseteq V_{i,0}$ such that all triples $\{v_a,v_{i+1},w\}$ with $1\le a \le j$ and $w\in V_{i,j}$ are in $\mathcal{E}$ or all are in $\mathcal{E}^c$. If the number of triples $\{v_{j+1},v_{i+1},w\}$ in $\mathcal{E}$ with $w\in V_{i,j}$ is at least $|V_{i,j}|/2$, then we set
$$
V_{i,j+1} = \{ w : \{v_{j+1},v_{i+1},w\} \in \mathcal{E} \text{ and } w\in V_{i,j}\},
$$
otherwise we set
$$
V_{i,j+1} = \{ w : \{v_{j+1},v_{i+1},w\} \in \mathcal{E}^c \text{ and } w\in V_{i,j}\}.
$$
Finally, we put $V_{i+1} = V_{i,i}$ and continue the algorithm until $A$ of size $m+1$ is chosen. This is possible, since $|V_i|\ge 1$ for all $1\le i\le m$. Indeed,
$$
|V_m| \ge \frac{|V_{m-1}|-1}{2^{m-1}} \ge \frac{|V_{m-1}|}{2^{m}} \ge \frac{|V_{m-2}|}{2^m 2^{m-1}}
\ge \frac{|V_1|}{2^m 2^{m-1}\cdots 2^2}
= \frac{n-1}{2^{(m+2)(m-1)/2}} \ge 1,
$$
since $n\ge 2^{m^2}$.

\section{Concluding remarks}
It seems natural to try to extend the above results for arbitrary $k\ge 4$. A similar approach to the one taken in Section~\ref{sec:lower} and inequality \eqref{eq:k} yield
\begin{equation}\label{eq:general}
c_1 (\underbrace{\log\log\dots\log}_{k-2} n)^{1/4}  \le f_{s,s+1}^{(k)} (n) \le c_2 (\log n)^{1/(k-2)}
\end{equation}
for some positive constants $c_1$ and $c_2$ depending only on $k$ and $s$. 
This big discrepancy could be possibly removed by strengthening the upper bound. Analogously to the problem of estimating the Ramsey numbers for hypergraphs one could apply some variation of the stepping-up lemma of Erd\H{o}s and Hajnal (see, e.g., Section~4.7 in~\cite{GRS}). Unfortunately, it is not obvious how to use this idea.

Recently, Conlon, Fox and Sudakov~\cite{CFS2} slightly improved the lower bound in \eqref{eq:general} and showed that 
$$
(\underbrace{\log\log\dots\log}_{k-2} n)^{1/3 - o(1)}  \le f_{s,s+1}^{(k)} (n). 
$$


\end{document}